\documentclass{amsart}

\usepackage{mathrsfs, amsfonts, amsmath, amssymb, amsthm}
\usepackage[cmtip,arrow]{xy}
\usepackage{pb-diagram,pb-xy}
\usepackage[dvips]{graphicx}
\usepackage{bbm}

\newtheorem{definition}{Definition}
\newtheorem{lemma}[definition]{Lemma}
\newtheorem{theorem}[definition]{Theorem}

\newtheorem{proposition}[definition]{Proposition}

\begin{document}

\title[Homological Stability of Non-orientable Mapping Class Groups]{Homological Stability of Non-orientable Mapping Class Groups with Marked Points}
\date{}
\author{Elizabeth Hanbury}

\address{Mathematical Institute, University of Oxford, 24-29 St Giles', Oxford, OX1 3LB}
\curraddr{Department of Mathematics, National University of Singapore, 2 Science Drive 2, Singapore 117543}
\email{liz.hanbury@gmail.com}


\begin{abstract}
Wahl recently proved that the homology of the non-orientable mapping class group stabilizes as the genus increases. In this short note we analyse the situation where the underlying non-orientable surfaces have marked points.
\end{abstract}

\maketitle

\section{Introduction}
Every non-orientable surface is homeomorphic to $S_{g,n}=\#_{g} \mathbb{R}P^{2} \setminus \coprod_{n}D^{2}$ for some unique $g,n$. The object of study in this note is the associated mapping class group $\mathcal{N}_{g,n}:=\pi_{0}\mathrm{Diff}(S_{g,n}, \partial)$ - the connected components of the group of diffeomorphisms of $S_{g,n}$ that fix the boundary pointwise. In \cite{Nathalie} Wahl proved that for $g \geq 4r +3$, $H_{r}(\mathcal{N}_{g,n})$ is independent of $g$ and $n$. This is the non-orientable version of Harer's well-known stability theorem for the orientable mapping class group \cite{Harer}.

Here we will look at non-orientable surfaces with marked points - the surface $S_{g,n}$ with $k$ marked points will be denoted $S_{g,n}^{k}$. We can consider the group of diffeomorphisms of $S_{g,n}^{k}$ that fix the marked points or the group of diffeomorphisms that may permute them. The corresponding mapping class groups will be denoted by $\mathcal{N}_{g,n}^{k}$ and $\mathcal{N}_{g,n}^{(k)}$ respectively. These groups fit into an exact sequence
\[
0 \to \mathcal{N}^{k}_{g,n} \to \mathcal{N}^{(k)}_{g,n} \to \Sigma_{k} \to 0.
\]

The first result of this note is the following corollary of Wahl's theorem.

\begin{theorem}
\label{oops}
For any $k \geq 0$, $H_{r}(\mathcal{N}_{g,n}^{k})$ and $H_{r}(\mathcal{N}_{g,n}^{(k)})$ are independent of $g$ and $n$ when $g$ is sufficiently large compared to $r$.
\end{theorem}

We show how to deduce this from Wahl's theorem. As she notes in her introduction, the original proof of stability could also be adapted to deal with marked points.

There are also stable versions of the non-orientable mapping class groups: for $n \geq 1$, $\mathcal{N}_{\infty,n}^{k}$ is defined to be the colimit of the system
\[
\mathcal{N}_{g,n}^{k} \stackrel{\alpha}{\to} \mathcal{N}_{g+1,n}^{k} \stackrel{\alpha}{\to} \mathcal{N}_{g+2,n}^{k} \stackrel{\alpha}{\to} \ldots
\]
where $\alpha$ is given by attaching $\mathbb{R}P^{2}\setminus (D^{2} \sqcup D^{2})$ to the underlying surface and extending diffeomorphisms by the identity. Similarly we have $\mathcal{N}_{\infty,n}^{(k)}$. Our second result is the following splitting theorem for the classifying spaces of these groups. This is analogous to Theorem 3.1 in \cite{Bod-Till} where the orientable case was considered.

\begin{theorem}
There are homology isomorphisms
\begin{eqnarray*}
B\mathcal{N}_{\infty,n}^{k} \to B\mathcal{N}_{\infty,1} \times B(O(2))^{k} & and &
B\mathcal{N}_{\infty,n}^{(k)} \to B\mathcal{N}_{\infty,1} \times B(\Sigma_{k} \wr O(2))
\end{eqnarray*}
\end{theorem}

In section \ref{suspect} we move on to look at how the homology of $\mathcal{N}_{g,n}^{(k)}$ depends on $k$. We prove
\begin{theorem}
\label{equal}
$H_{r}(\mathcal{N}_{g,n}^{(k)})$ is independent of $k$ (as well as $g$ and $n$) when $g$ and $k$ are sufficiently large relative to $r$.
\end{theorem}
Our techniques could also be applied to the orientable case to prove
\begin{theorem}
\label{more-equal}
$H_{r}(\Gamma_{g,n}^{(k)})$ is independent of $k$ when $g$ and $k$ are sufficiently large relative to $r$.
\end{theorem}
Here $\Gamma_{g,n}^{(k)}$ is the orientable mapping class group in which diffeomorphisms are allowed to permute marked points.

Hatcher and Wahl have recently proved versions of Theorems \ref{equal} and \ref{more-equal} in which they don't require the genus to be large \cite{HW}. They obtain these results as special cases of the more general theory of homological stability of mapping class groups of $3$-manifolds.

The proof of Theorem \ref{equal} is based on the following theorem on the Borel construction. This theorem is well-known but we are not aware of any published proof.
\begin{theorem}
\label{buttmunch}
If $X$ is a path-connected space then for $k \geq 2r+6$ there is an isomorphism
\[
H_{r}(E\Sigma_{k} \times_{\Sigma_{k}} X^{k}) \to H_{r}(E\Sigma_{k+1} \times_{\Sigma_{k+1}} X^{k+1}).
\]
\end{theorem}
This theorem is based on the homological stability of symmetric groups with twisted coefficients \cite{Betley}. The proof is contained in the appendix.

\vspace{\headheight}
\section{Homological Stability with Marked Points}

For any $n\geq 1$ there is a stabilization map $\alpha:\mathcal{N}_{g,n} \to \mathcal{N}_{g+1,n}$ as described in the introduction. Similarly there are maps $\beta:\mathcal{N}_{g,n} \to \mathcal{N}_{g,n+1}$, given by attaching a pair of pants, and $\gamma:\mathcal{N}_{g,n} \to \mathcal{N}_{g,n-1}$ given by attaching a disc to one of the boundary components. The map $\gamma$ forms a left-inverse to $\beta$ for $n \geq 2$. The main result of \cite{Nathalie} is the following.

\begin{theorem}[Wahl]
\label{julie}
For any $n\geq 1$
\begin{itemize}
\item[(i)] $\alpha_{\ast}: H_{r}(\mathcal{N}_{g,n}) \to H_{r}(\mathcal{N}_{g+1,n})$ is an isomorphism provided $g \geq 4r+3$
\item[(ii)] $\beta_{\ast}: H_{r}(\mathcal{N}_{g,n}) \to H_{r}(\mathcal{N}_{g,n+1})$ is an isomorphism provided $g \geq 4r+3$
\item[(iii)] $\gamma_{\ast}: H_{r}(\mathcal{N}_{g,1}) \to H_{r}(\mathcal{N}_{g,0})$ is an isomorphism provided $g \geq 4r+5$.
\end{itemize}
\end{theorem}

The stabilization maps $\alpha, \beta$ and $\gamma$ also define maps on $\mathcal{N}_{g,n}^{k}$ and $\mathcal{N}_{g,n}^{(k)}$. Theorem \ref{oops} follows from the following extension of Wahl's theorem.

\begin{theorem}
\label{eat}
Theorem \ref{julie} still holds if we replace $\mathcal{N}_{g,n}$ by $\mathcal{N}_{g,n}^{k}$ or $\mathcal{N}_{g,n}^{(k)}$.
\end{theorem}

\begin{proof}
We will give the proof of (i) for $\mathcal{N}_{g,n}^{(k)}$ - the other cases are entirely analogous.

Let $\mathrm{Diff} (S_{g,n}^{(k)}, \partial)$ denote the group of diffeomorphisms of $S_{g,n}^{k}$ which fix the boundary and may permute the marked points. Fix a framing at each marked point i.e. a pair of orthonormal tangent vectors. By considering the effect of a diffeomorphism on the marked points and their framings we obtain a map
\begin{equation}
\label{squidge}
\mathrm{Diff} (S_{g,n}^{(k)},\partial) \stackrel{\rho}{\to} \Sigma_{k} \wr Gl_{2} (\mathbb{R}) \simeq \Sigma_{k} \wr O(2)
\end{equation}
which is a fibration onto its image. The fiber is homotopic to the group of diffeomorphisms which fix a neighbourhood of each marked point i.e. $\mathrm{Diff}(S_{g,n+k},\partial)$.

We claim that $\rho$ is in fact surjective. By considering diffeomorphisms of $D^{2}$ with $k$ points removed we can see that every element of $\Sigma_{k} \wr Gl_{2}^{+}(\mathbb{R})$ lies in the image - we complete the proof of the claim by showing that $\iota \times - \mathrm{id}$ also lies in the image where $\iota$ denotes the identity permutation. This amounts to displaying a diffeomorphism of $S_{g,n}^{k}$ which reverses the orientation at a chosen marked point and fixes everything else. Such a diffeomorphism is given by a `dragging function': fix a curve $C$ which passes through the chosen marked point and such that a neighbourhood $U$ of $C$ is homeomorphic to a M\"obius band. Choose $U$ to avoid the other marked points. The diffeomorphism we require is obtained by dragging the marked point all the way around $C$ and is the identity outside $U$. Explicitly we can define this as a composition of $m$ diffeomorphisms $\{ \epsilon_{j}\}_{j=1}^{m}$ which each move the marked point $\frac{1}{m}$th of the way around $C$ and are the identity outside $U$. See Figure 1 for a representation of $\epsilon_{j}$.

\begin{figure}
\includegraphics{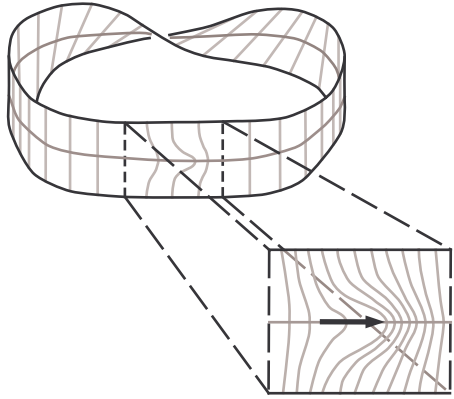}
\caption{The diffeomorphism $\epsilon_{j}$}
\label{}
\end{figure}

Now we take classifying spaces in (\ref{squidge}) and use the fact that the components of the diffeomorphism groups are contractible \cite{EE, ES} to obtain a homotopy fibration
\begin{equation}
\label{excess}
B\mathcal{N}_{g,n+k} \to B\mathcal{N}_{g,n}^{(k)} \to B(\Sigma_{k} \wr O(2)).
\end{equation}
The map $\alpha$ induces a map from this fibration to
\[
B\mathcal{N}_{g+1,n+k} \to B\mathcal{N}_{g+1,n}^{(k)} \to B(\Sigma_{k} \wr O(2))
\]
which is the identity on the base spaces.

The induced map on the $E^{2}$-page of the associated Serre spectral sequences is the obvious map
\begin{eqnarray*}
E_{p,q}^{2} = H_{p}(B(\Sigma_{k} \wr O(2)); H_{q}(B\mathcal{N}_{g,n+k})) \\
\to \tilde{E}_{p,q}^{2} = H_{p}(B(\Sigma_{k} \wr O(2)); H_{q}(B\mathcal{N}_{g+1,n+k})).
\end{eqnarray*}

By Theorem \ref{julie} the map on the coefficients is an isomorphism for $g \geq 4q+3$. Hence the map on the $E_{2}$-term is an isomorphism for all $p$ and for $q \leq \frac{g-3}{4}$ so by the Zeeman comparison theorem (see section 1 in \cite{Ivanov}) the map on the abutments
\[
H_{r}(\mathcal{N}_{g,n}^{(k)}) \to H_{r}(\mathcal{N}_{g+1,n}^{(k)})
\]
is an isomorphism for $r \leq \frac{g-3}{4}$.
\end{proof}

\vspace{\headheight}
\section{Splitting of the Classifying Spaces}

In this section we will establish a splitting result for the classifying spaces of $\mathcal{N}_{\infty,n}^{k}$ and $\mathcal{N}_{\infty,n}^{(k)}$ cf. \cite{Bod-Till}.
\begin{theorem}
\label{patricia}
There are homology isomorphisms
\begin{itemize}
\item[(i)] $B\mathcal{N}_{\infty,n}^{k} \to B\mathcal{N}_{\infty,1} \times B(O(2))^{k}$
\item[(ii)] $B\mathcal{N}_{\infty,n}^{(k)} \to B\mathcal{N}_{\infty,1} \times B(\Sigma_{k} \wr O(2))$.
\end{itemize}
\end{theorem}

\begin{proof}
We give the proof of (ii), the proof of (i) is completely analogous. Consider the fibration (\ref{excess}) from the previous section
\[
B\mathcal{N}_{g,n+k} \to B\mathcal{N}_{g,n}^{(k)} \stackrel{B\rho}{\longrightarrow} B(\Sigma_{k} \wr O(2)).
\]
There is a map from this fibration to the trivial fibration
\[
B\mathcal{N}_{g,1} \to B\mathcal{N}_{g,1} \times B(\Sigma_{k} \wr O(2)) \to B(\Sigma_{k} \wr O(2)).
\]
On the total space this is given in the first factor by forgetting marked points and applying $\gamma$ and in the second factor by $B\rho$. On the fiber it's induced by $\gamma$.

We know from Theorem \ref{julie} that the maps on the fiber and on the base are $H_{r}$-isomorphisms for $g$ sufficiently large. Considering the associated Serre spectral sequences and applying the Zeeman comparison theorem we are able to conclude that the map on total spaces is an $H_{r}$-isomorphism for sufficiently large $g$. Letting $g$ tend to infinity we obtain the result.

\end{proof}

\section{Changing the Number of Marked Points}
\label{suspect}

We now move on to investigating how the homology of $\mathcal{N}_{g,n}^{k}$ and $\mathcal{N}_{g,n}^{(k)}$ depend on $k$. We can see from Theorem \ref{patricia} that the homology of $\mathcal{N}_{g,n}^{k}$ is dependent on $k$ no matter how large $g$ or $k$ become. However we have maps $\mu:\mathcal{N}_{g,n}^{k} \to \mathcal{N}_{g,n}^{k+1}$ and $\nu:\mathcal{N}_{g,n}^{k+1} \to \mathcal{N}_{g,n}^{k}$ where $\mu$ is given by gluing on a cylinder with a marked point and extending maps by the identity and $\nu$ is given by forgetting the $(k+1)$st marked point. The map $\nu$ is a left-inverse to $\mu$ so we have

\begin{lemma}
$\mu_{\ast}:H_{r}(\mathcal{N}_{g,n}^{k}) \to H_{r}(\mathcal{N}_{g,n}^{k+1})$ is injective.
\end{lemma}

When we are dealing with $\mathcal{N}_{g,n}^{(k)}$ we can say much more. In fact in the next theorem we will prove that $H_{r}(\mathcal{N}_{g,n}^{(k)})$ is independent of $k$ provided $g$ and $k$ are sufficiently large. The proof relies on Betley's result on the homological stability of the symmetric groups with twisted coefficients \cite{Betley}.

\begin{theorem}
\label{cherie}
When $g \geq 4r+3$ and $k \geq 2r+8$, gluing on a cylinder with a marked point induces an isomorphism
\[
H_{r}(\mathcal{N}_{g,n}^{(k)}) \stackrel{\cong}{\longrightarrow} H_{r}(\mathcal{N}_{g,n}^{(k+1)}).
\]
\end{theorem}

\begin{proof}
Recall the fibration (\ref{excess}) from above
\[
B\mathcal{N}_{g,n+k} \to B\mathcal{N}_{g,n}^{(k)} \longrightarrow B(\Sigma_{k} \wr O(2)).
\]
Gluing on a cylinder with a marked point induces a map from this fibration to
\[
B\mathcal{N}_{g,n+k+1} \to B\mathcal{N}_{g,n}^{(k+1)} \longrightarrow B(\Sigma_{k+1} \wr O(2)).
\]
There is an induced map of the associated Serre spectral sequences and on the $E^{2}$-term  it is
\begin{eqnarray*}
\label{criticise}
E^{2}_{p,q} = H_{p}(B(\Sigma_{k} \wr O(2)); H_{q}(\mathcal{N}_{g,n+k})) \\
\to \tilde{E}^{2}_{p,q} = H_{p}(B(\Sigma_{k+1} \wr O(2)); H_{q}(\mathcal{N}_{g,n+k+1})).
\end{eqnarray*}
The map on the coefficients is an isomorphism for $g \geq 4q+3$ by Theorem \ref{julie}. We also have that the $\pi_{1}(B(\Sigma_{k} \wr O(2)))$-action on $H_{q}(\mathcal{N}_{g,n+k})$ is trivial for $g$ in this range - this is proved in the same way as Lemma 3.3 in \cite{Bod-Till} - so the coefficients are not twisted. Since $B(\Sigma_{k} \wr O(2)) \simeq E\Sigma_{k} \times_{\Sigma_{k}} BO(2)^{k}$ then, for $g \geq 4q+3$, we can apply Theorem \ref{buttmunch} and deduce that the map on the $E^{2}$-term is an isomorphism for $p \leq \frac{k-6}{2}$, $q \leq \frac{g-3}{4}$. Applying the Zeeman comparison theorem yields the result.
\end{proof}

\vspace{1mm}
\section{Appendix: Applying Betley's Stability Result}

Let $\Phi$ be the category whose objects are the finite pointed sets $[n]=\{ 0,1,\ldots, n \}$ and whose morphisms are pointed maps with the property that the pre-image of each non-zero element contains at most one element. Let $Ab$ denote the category of abelian groups. A coefficient system is a functor $T: \Phi \to Ab$ with $T([0])=0$. Thus $\Sigma_{n}$ acts on $T([n])$ for each $n$ and we have the homology groups $H_{i}(\Sigma_{n}; T([n]))$.

In order to define the notion of degree of coefficient systems, we
consider the maps $e_{s}:[n] \to [n-1]$ for $1 \leq s \leq n$ given
by $e_{s}(j)=j$ for $j < s$, $e_{s}(s)=0$ and $e_{s}(j)=j-1$ for
$j>s$.  The coefficient system $T$ is said to be of degree $d$ if
\[
T_{d+1}:= ker \left( \prod_{s=1}^{d+1}(e_{s})_{\ast}: T([d+1]) \to \prod_{s=1}^{d+1}T([d]) \right)
\]
is zero but $T_{d}$ is non-zero.

Betley's main result was the following
\begin{theorem}[Betley]
\label{spam}
For any coefficient system $T$ of degree $\leq d$ the natural map
\[
H_{i}(\Sigma_{n}; T([n])) \to H_{i}(\Sigma_{n+1}; T([n+1]))
\]
is an isomorphism for $2i + d \leq n$.
\end{theorem}

Betley actually defined a coefficient system to be a functor on the category $\Gamma$ of finite pointed sets with pointed maps, which contains $\Phi$. Thus Theorem \ref{spam} is a slightly stronger result than Betley originally proved. We'll summarize his proof to check that it's still valid in this more general situation.

\begin{proof}[Proof of Theorem \ref{spam}.]
There is an isomorphism of abelian groups;
\begin{equation}
\label{angel}
T([n]) \cong \bigoplus_{k} T^{k}_{n}
\end{equation}
where for given $k$, $T^{k}_{n}=\bigoplus_{\gamma}T_{k}$, the latter sum being taken over all order-preserving, pointed injections $\gamma:[k] \to [n]$. This is formula (1.1.1) in Betley's paper. A proof can be found in \cite{Eilenberg-MacLane}, Theorem 9.1, and only requires the morphisms from $\Phi$.

The next step is to identify the $\Sigma_{n}$-action on $\bigoplus_{k}T^{k}_{n}$ under (\ref{angel}). The action respects the decomposition. Each $T_{k} \leq T([k])$ has a natural $\Sigma_{k}$-action. Careful consideration of the bijection (\ref{angel}) gives an isomorphism of $\Sigma_{n}$-modules
\begin{equation}
\label{arithmetic}
T_{n}^{k} \cong T_{k} \otimes_{\mathbb{Z}[\Sigma_{k} \times \Sigma_{n-k}]} \mathbb{Z}[\Sigma_{n}]
\end{equation}
where $\Sigma_{n-k}$ acts trivially on $T_{k}$. This is in the proof of Betley's Lemma 3.1.

Now if $T$ has degree $d$ and $n>d$ we have
\[
H_{i}(\Sigma_{n};T([n])) \cong H_{i}(\Sigma_{n}; \bigoplus_{k=1}^{d} T_{n}^{k}) \cong \bigoplus_{k=1}^{d} H_{i}(\Sigma_{n}; T_{n}^{k}).
\]
By (\ref{arithmetic}), Shapiro's lemma and the K\"unneth theorem, this is isomorphic to
\[
\bigoplus_{k=1}^{d} H_{i}(\Sigma_{k} \times \Sigma_{n-k}; T_{k}) \cong \bigoplus_{k=1}^{d} \{ \bigoplus_{j+l=i} H_{j} (\Sigma_{n-k}; H_{l}(\Sigma_{k};T_{k})) \}.
\]

The same formula holds for $H_{i}(\Sigma_{n+1}; T([n+1]))$. The stabilization map is induced by the canonical inclusion $\Sigma_{n-k} \to \Sigma_{n-k+1}$ and is the identity on the coefficients. The result is now immediate by the homological stability of the symmetric groups.
\end{proof}

The proof of Theorem \ref{buttmunch} requires the following proposition
\begin{proposition}
\label{balula}
If $X$ is a path-connected space with basepoint and
$F$ is any field then $[n] \mapsto H_{q}(X^{n},F)$ defines a
coefficient system of degree $\leq q$.
\end{proposition}

\begin{proof}
For ease of notation, we'll suppress the coefficient field $F$. First we'll check that this assignment has the necessary functoriality. Let $\lambda:[k] \to [n]$ be a morphism in $\Phi$, there is an induced map $\lambda_{\ast}:X^{k} \to X^{n}$ given by $(x_{1}, \ldots, x_{k}) \to (x_{\lambda^{-1}(1)}, \ldots, x_{\lambda^{-1}(n)})$ where $x_{\lambda^{-1}(j)}:=\ast$ if $\lambda^{-1}(j)=\emptyset$. Clearly $(\lambda \circ \mu)_{\ast}=\lambda_{\ast} \circ \mu_{\ast}$ so this gives us functoriality.

Now we move on to checking that the coefficient system has finite degree. By the definition of degree, our task is to prove that the kernel of the map
\begin{equation}
\label{joseph}
\prod_{s=1}^{n}(e_{s})_{\ast}: H_{q}(X^{n}) \to \prod_{s=1}^{n} H_{q}(X^{n-1})
\end{equation}
is trivial when $n > q$.

Since we're working with field coefficients, the K\"unneth theorem tells us that,
\[
H_{q}(X^{n}) \cong \bigoplus_{i_{1}+ \ldots + i_{n}=q} H_{i_{1}}(X) \otimes_{F} \ldots \otimes_{F} H_{i_{n}}(X).
\]

The induced map $(e_{s})_{\ast}:H_{q}(X^{n}) \to  H_{q}(X^{n-1})$
takes the summand corresponding to $(i_{1}, \ldots, i_{n})$ to zero
if $i_{s} \neq 0$ and isomorphically to the summand corresponding to
$(i_{1}, \ldots, \hat{i_{s}}, \ldots, i_{n})$ if $i_{s}=0$. To see
this, think of $(e_{s})_{\ast}: X^{n} \to X^{s-1} \times \{\ast\}
\times X^{n-s} \cong X^{n-1}$ as the map sending the $s$th factor to
the basepoint. Then the induced map on homology takes $H_{i_{r}}(X)$
by the identity to $H_{i_{r}}(X)$ when $r \neq s$ and takes
$H_{i_{s}}(X)$ to $H_{i_{s}}(\ast)$ by the obvious induced map. The
latter is an isomorphism if $i_{s}=0$ (because $X$ is
path-connected) and is the zero map otherwise. Thus the kernel of
$(e_{s})_{\ast}$ is
\[
\bigoplus_{\stackrel{i_{1}+ \ldots + i_{n}=q}{i_{s} \neq 0}} H_{i_{1}}(X) \otimes_{F} \ldots \otimes_{F} H_{i_{n}}(X)
\]
The kernel of $\prod_{s=1}^{n}(e_{s})_{\ast}$ is the intersection of these kernels and hence it is the direct sum of those summands corresponding to $n$-tuples $(i_{1}, \ldots, i_{n})$ with $i_{1}+ \ldots + i_{n}=q$ and every $i_{s} \neq 0$. When $n > q$ there are no such $n$-tuples so the kernel is trivial.
\end{proof}

\vspace{2mm}
\begin{proof}[Proof of Theorem \ref{buttmunch}.]
Consider the map of fibrations
\begin{eqnarray*}
\begin{diagram}
\node{X^{k}} \arrow{e} \arrow{s} \node{E\Sigma_{k} \times_{\Sigma_{k}} X^{k}} \arrow{e} \arrow{s} \node{B\Sigma_{k}} \arrow{s}\\
\node{X^{k+1}} \arrow{e} \node{E\Sigma_{k+1} \times_{\Sigma_{k+1}} X^{k+1}} \arrow{e} \node{B\Sigma_{k+1}}
\end{diagram}
\end{eqnarray*}
and the associated Serre spectral sequences with coefficients in some field $F$. There is an induced map of the spectral sequences and on the $E^{2}$-term it is given by
\[
E_{p,q}^{2} = H_{p}(\Sigma_{k}; H_{q}(X^{k},F)) \to \tilde{E}_{p,q}^{2} = H_{p}(\Sigma_{k+1}; H_{q}(X^{k+1},F)).
\]
But by Proposition \ref{balula}, $\{H_{q}(X^{k},F)\}_{k}$ is a coefficient system of degree $\leq q$. Hence Betley's stability result implies that this map is an isomorphism for $2p+q \leq k$. Now by the Zeeman comparison theorem, the induced map on the abutments
\[
H_{r}(E\Sigma_{k} \times_{\Sigma_{k}} X^{k},F) \to H_{r}(E\Sigma_{k+1} \times_{\Sigma_{k+1}} X^{k+1},F)
\]
is an isomorphism for $r \leq \frac{k}{2}-1$. Since $E\Sigma_{k} \times_{\Sigma_{k}} X^{k} \to E\Sigma_{k+1} \times_{\Sigma_{k+1}} X^{k+1}$ induces an isomorphism on $H_{r}(-,F)$ for any field $F$ when $r \leq \frac{k}{2}-1$ it also induces an isomorphism on $H_{r}(-,\mathbb{Z})$ for $r \leq \frac{k}{2}-3$. This is a consequence of the universal coefficient theorem.
\end{proof}

\bibliographystyle{plain}
\bibliography{thesis}

\end{document}